%% file: article_multigrid_dwarka_2023.tex
\documentclass[onefignum,onetabnum]{siamart190516}
\usepackage{mathtools}
\usepackage{caption}
\usepackage{subfigure}
\usepackage{amsmath,amsfonts,epsfig}
\usepackage{tikz}
\usetikzlibrary{positioning}
\usetikzlibrary{matrix}

\usepackage{tikz-cd}
\usepackage[all,cmtip]{xy}
\usepackage{mathabx} 
\usepackage{dsfont}
\usepackage{graphicx}
\usepackage{cleveref}
\usepackage{mwe}
\usepackage{graphicx}
\usepackage{booktabs}
\usepackage[utf8]{inputenc}
\usepackage{color}
\usepackage{placeins}

\usepackage{sidecap}
\usepackage[colorinlistoftodos]{todonotes}
\reversemarginpar
\input{ex_shared}

\ifpdf
\hypersetup{
  pdftitle={Stand-alone multigrid for Helmholtz revisited: towards convergence using standard components},
  pdfauthor={Vandana Dwarka, Cornelis Vuik}
}
\fi

\begin{document}

\maketitle

% REQUIRED
\begin{abstract}
Getting standard multigrid to work efficiently for the high-frequency Helmholtz equation has been an open problem in applied mathematics for years. Much effort has been dedicated to finding solution methods which can use multigrid components to obtain solvers with a linear time complexity. 
In this work we present one among the first stand-alone multigrid solvers for the 2D Helmholtz equation using both a constant and non-constant wavenumber model problem. We use standard smoothing techniques and do not impose any restrictions on the number of grid points per wavelength on the coarse-grid. As a result we are able to obtain a full V- and W-cycle algorithm.
The key features of the algorithm are the use of higher-order inter-grid transfer operators combined with a complex constant in the coarsening process. Using weighted-Jacobi smoothing, we obtain a solver which is $h-$independent and scales linearly with the wavenumber $k$. Numerical results using 1 to 5 GMRES(3) smoothing steps approach $k-$ and $h-$ independent convergence, when combined with the higher-order inter-grid transfer operators and a small or even zero complex shift. The proposed algorithm provides an important step towards the perpetuating branch of research in finding scalable solvers for challenging wave propagation problems. 
\end{abstract}

% REQUIRED

% REQUIRED
\begin{AMS}
  {\color{black} 65N55, 65F10, 65F15}
\end{AMS}
%\vspace*{-0.25cm}
%%%%%%%%%%%%%%%%%%%%%%%  START SECTION 1 %%%%%%%%%%%%%%%%%%%%%%%

%%%%%%%%%%%%%%%%%%%%%%%%%%%%%%%%%%%%%%%%%%%%%%%%
%%%%%%%%%%%%%%%%%%%%%%%%%%%%%%%%%%%%%%%%%%%%%%%%
%\newpage

\noindent 

\section{Introduction} %\label{sec1}  %% 11-5
For over  three decades, researchers have studied ways to design convergent and efficient multigrid schemes for the indefinite Helmholtz equation, which to this day remains an open problem in applied mathematics. The Helmholtz equation finds its application in many engineering practices, ranging from geoscience to plasma physics and quantum mechanics. Solving the Helmholtz equation numerically is a challenging task. Multigrid methods are widely known to be efficient and scalable numerical solvers for elliptic PDEs. It provided a state-of-the-art numerical approach, where the amount of computational work scales with the number of unknowns. In essence, the method combines solutions on coarser grids with relaxation (smoothing) techniques to arrive at a solution on the fine grid. 

While multigrid works well for the positive-definite variant of the Helmholtz equation or when the wavenumber is small, the method generally fails to converge as a stand-alone solver for the indefinite version. This break-down can be attributed to two things. First of all, the coarse-grid needs to be fine enough in order to resolve the waves. A difference between the wavenumber on the fine vs. coarse-grid(s) can cause a phase error leading to instabilities \cite{ernst2012difficult}. Moreover, eigenvalues approaching zero also hamper the convergence as they lead to near-singularity within the coarse-grid operator \cite{ernst2012difficult}. 
Several strategies have been proposed to study the convergence behavior of such indefinite systems \cite{brandt1986multigrid, bramble1988analysis, wang1993convergence, hackbusch2013multi}. This boils down to prove convergence of the two-grid method for the positive definite variant (see \cite{hackbusch2013multi}, chapter 11.6). Another approach lies in trying to find upper and lower bounds to the field of values of the two-grid operator \cite{notay2010algebraic,notay2007convergence,notay2010aggregation,notay2020analysis}. This requires that the Hermitian part of the linear system $\frac{1}{2}\left(A + A^{T} \right)$ is positive definite, which is not the case for the indefinite Helmholtz equation. While we can still use the condition that the spectral radius of the two-grid iteration matrix should be smaller than one, the computation of the eigenvalues becomes difficult in case of Sommerfeld boundary conditions. {\color{black}Moreover, it has been noted that for nonnormal (complex) systems in particular, the Helmholtz equation being an example of such, assessing convergence solely based on the spectral radii can be misleading as it only rigorously describes asymptotic convergence \cite{notay2020analysis}. Bounds in terms of norms would be more appropriate. However, these often go hand in hand with assumptions which are difficult to check and contain unknown constants, which reduces their practical applicability.}
 
For positive-definite and complex Hermitian problems, Geometric and/or Algebraic Multigrid are powerful solvers which provide solutions in linear time complexity by constructing coarser linear systems based on the information provided by the fine-grid linear system. For complex Hermitian matrices, the coarse linear systems can be build using the Galerkin condition, which has been studied in \cite{lahaye2000algebraic, maclachlan2008algebraic}. If the use of a simple smoother by means of a stationary iterative method such as weighted-Jacobi or Gauss-Seidel is preferred, the coarse grid correction scheme should balance the simplicity of the smoother by accurately transferring smooth errors from the coarse grid to the fine grid \cite{maclachlan2008algebraic}. When the linear system is complex symmetric, linear interpolation can be used to construct the prolongation operator. The restriction operator can then be chosen as the transpose of the prolongation operator. In fact, for real matrices optimality studies have been conducted in order to construct the interpolation and prolongation operators \cite{brezina2001algebraic, brezina2005adaptive, falgout2004generalizing,brannick2018optimal}. In \cite{ramos2019optimal}, a similar optimality condition from \cite{falgout2004generalizing,brannick2018optimal} has been extended to the Hermitian case. These studies show that the optimal prolongation and restriction operator are based on the eigenvectors of the system matrix and several approximations are given which could function as a satisfactory proxy in the case of symmetric positive definite (SPD) and Hermitian systems respectively. 

Including a complex shift, leads to the Complex Shifted Laplacian (CSL), which can be solved using multigrid as a stand-alone solver. Consequently, for the original indefinite Helmholtz problem, multigrid applied to the CSL can be used as a preconditioner, {\color{black} together with an iterative Krylov subspace method, such as {\color{black}GMRES}}. Here, the original indefinite Helmholtz problem is solved by approximately inverting the CSL using geometric multigrid \cite{erlangga2006novel,erlangga2008multilevel,erlangga2012iterative, cocquet2017large}. While this works well for homogeneous problems, it is not suitable for heterogeneous problems with high contrasts \cite{sheikh2016accelerating,erlangga2008multilevel}. 

When it comes to using multigrid as a stand-alone solver for the original Helmholtz operator, several studies have focused on getting multigrid to comply with the underlying characteristics of the indefinite nonsymmetric system. 
One study uses a smoothed aggregation approach, where an effective prolongation operator is constructed by minimizing the energy norm defined using the normal equations $A^{*}A$. This is achieved by running the CGNR algorithm to solve the equation $A^{*}AP = 0$, where $P$ is the tentative prolongation operator \cite{olson2010smoothed,olson2011general}. Two other methods are the wave-ray approach \cite{livshits2006accuracy} and the first-order least squares (FOLS) approach \cite{lee2000first}. Both compute multiple coarse spaces to approximate the plane waves in the near null space. In both works, the development of more efficient approaches for such linear systems has been stressed. Moreover, for variable wave numbers, it remains challenging to choose an effective prolongation and restriction operator, especially if the aim is to keep the algorithm as light as possible. 

Another difficulty arises in picking a suitable smoother, as for Helmholtz problems, the (weighted)-Jacobi method diverges. While the latter has natural parallelization properties and is often preferred for SPD problems, no clear remedy exists to apply it suitably for Helmholtz problems at large wavenumbers and tailored adaptions are needed \cite{hadley2005complex,ernst2013multigrid,erlangga2006novel}. Essentially, this relates to the indefiniteness of the underlying operators, leading to negative eigenvalues which ensure that the spectral radius of the weighted-Jacobi iteration matrix is always greater than 1. One important recent study finds that replacing the weighted-Jacobi smoother with a few {\color{black}GMRES} iterations could lead to improved convergence \cite{elman2001multigrid}. In \cite{ernst2012difficult,ernst2013multigrid}, it was shown that for the one-dimensional indefinite Helmholtz equation, convergence can be achieved by considering a two-step Jacobi smoother. Scheduling the smoothing steps on all levels as a function of the wave number provides a uniform error reduction in the high-frequency modes. More recently a two-dimensional dispersion correction has been developed for the original Helmholtz operator when using a constant wave number. Here, a larger number of smoothing steps combined with at least 4.5 grid points per wavelength at the coarsest level leads to a convergent multigrid scheme. \cite{cocquet2021closed}.

A full multigrid hierarchy, without any restrictions on the grid resolution at the coarsest level, is given in \cite{cools2014new}. Here, a re-discretization technique is used where a level-dependent complex shift is incorporated, except the finest level. The authors show that the altered two-grid correction scheme, combined with 3 GMRES iterations as a smoother, remains robust and does not lead to agitation of smooth modes. The method works well on structured grids and in cases where the behavior of the wave number on the fine and coarse(r) grids can be defined a-priori. {\color{black}Complex geometries, high-contrast and irregularly varying wave number would pose difficulties in constructing the coarser linear systems along the multilevel hierarchy, when re-discretization is used combined with a level-dependent parameter.}   

In this work we aim to solve the complex indefinite Helmholtz problem using constant and variable wave numbers with a stand-alone multigrid method using Galerkin coarsening. {\color{black}For the baseline convergence study, we impose \textit{two} requirements: the method should converge with a small number of smoothing steps and the intergrid transfer operators should be level-independent and easy to construct and implement. 
% \begin{enumerate}
%     \item The method should converge with a small number of smoothing steps
%     \item The prolongation and restriction operators should be easy to construct, implement and be level independent.
% \end{enumerate}
}
\newline
{\color{black}
Our contribution to this long standing field of research can be summarized into the following points:
\begin{itemize}
    \item To the best of our knowledge, we are the \textit{first} to present a fully convergent multigrid solver for the two-dimensional indefinite nonsymmetric Helmholtz equation using classical multigrid components for elliptic problems (V- and W-cycle hierarchy, damped Jacobi smoothing and no restriction on the (coarsest) levels.
    \item We provide a novel generalized proof to assess two-grid convergence for any (highly) indefinite system. Our proof is not limited to the Helmholtz equation and does not require the calculation of eigenvalues. 
     \item The novelty in our theory also lies in the fact that we do not require assumptions on normality of the operator. 
    \item For Helmholtz problems in particular, we provide a optimality condition which can easily be verified without computing eigenvalues. 
    \item The results are robust to large wavenumbers and non-constant wavenumbers representing high-contrast heterogeneities. 
\end{itemize}

% \begin{itemize}
% \item proof general complex indefinite non-adjoint/Hermitian systems, not necessarily only for Helmholtz
% \item Hackbush connectie beter duidelijk maken
% \item Optimality condition
% \item Emphasis dat het gaat niet om hoe snel het convergeert, als het maar convergeert. De number of iterations is dan van latere zorg. 
% \end{itemize}

In particular, we show that the use of coarsening on the CSL combined with a higher-order prolongation and restriction operator leads to a convergent scheme, even when using the classical divergent weighted-Jacobi smoother. 

Several remarks are in place. First of all, in this work we solely consider getting a convergent solver rather than a solver with wavenumber independent convergence. This is a topic for future research due to the sheer difficulty in getting a convergent stand-alone multigrid solver in the first place. Second of all, given that multigrid generally already diverges for one-dimensional Helmholtz problems with small wavenumbers, this research will, in accordance with the literature, focus on one- and two-dimensional model problems. The study of three-dimensional problems is a topic for future research as it requires a thorough understanding of the convergence behavior in one- and two-dimensions first. 

% benadrukken al in 1d en 2d conv bereiken onmogelijk -> future research 3d. 
% The main new contribution of this paper is that no restriction is imposed on the number of grid points per wavelength at the coarsest level and no level-dependent parameters are introduced, keeping the general configuration similar to the original standard multigrid algorithm for elliptic problems. As a result, we are able to construct a full V- and W-cycle up to reaching only one grid point. Moreover, the theory developed in this work is novel in the sense that it allows for the assessment of convergence for highly indefinite systems without using any assumptions as regards adjointness or convergent smoothers.
}

{\color{black}The structure of this paper is as follows. In \cref{sec2} we introduce the model problems which will be the focus of this work. In \cref{sec13} we give an introduction to the multigrid scheme and its parameters, to continue with the analysis of the convergence of multigrid methods for Helmholtz problems in \cref{sec3}. Numerical experiments are presented in \cref{sec5}. Our conclusions are summarized in \cref{sec6}.} 

\section{Problem Description} \label{sec2}
We continue by defining the model problems which will be studied in this paper. We focus on the 2D constant wavenumber model using $k > 0$, which we call MP 2-A. 
\begin{eqnarray*} %\label{MP1A}
-\Delta u(x,y) - k^2u(x,y) \, u(x,y) &=&  \delta{(x - \frac{1}{2},y - \frac{1}{2})}, \quad {(x,y)} \in  \Omega = [0,1]^2,  \\
\left( \frac{\partial{}}{\partial{{\mathbf{n}}}} - ik \right)  u({x,y}) &=& 0, \quad (x,y) \in \partial{\Omega}. 
\end{eqnarray*}

The final test problem uses a non-constant wavenumber $k(x,y)$. This gives us MP 2-B as defined below
%\subsubsection{MP 2-B}
\begin{eqnarray*} %\label{MP1A}
-\Delta u(x,y) - {k(x,y)}^2u(x,y) \, u(x,y) &=&  \delta{(x - \frac{1}{2},y - \frac{1}{2})}, \quad {(x,y)} \in  \Omega = [0,1]^2,  \\ 
\left( \frac{\partial{}}{\partial{{\mathbf{n}}}} - ik(x,y) \right)  u({x,y}) &=& 0, \quad (x,y) \in \partial{\Omega},  
\end{eqnarray*}
where $\chi(x,y)$ is a random real function in the range of $[0,1]$ and $k_1,k_2$ are positive real numbers. Note that in this case $k$ varies exactly between $k_1$ and $k_2$. We use second-order finite differences to discretize the model problems and define the step-size $h = \frac{1}{n}$, where $n$ is chosen such that for each $k_1,k_2$ we have $kh = 0.625$. This is equivalent to having 10 grid points per wavelength.

\section{Multigrid Methods} \label{sec13}  %% 11-5
Multigrid method as a stand-alone solver combines two ingredients in a recursive way: a coarse grid correction and a smoother. If we let $P$ and $P'$ denote the matrix variant of the prolongation and restriction operators, then we define the coarse grid correction matrix $CGC$ as 
\begin{align*}
    CGC = I - P'{A_c}^{-1}PA, \quad A_c = PAP'.
\end{align*}
Note that we here denote the coarse grid version of $A$ by $A_c$ instead of $A_H$ as the matrix $A$ is now complex due to the Sommerfeld boundary conditions and we reserve the $H-$script to denote the conjugate transpose. 

As for the smoother, we work with the standard weighted-Jacobi smoother, which is defined as
\begin{align*}
    S_{JAC} = I - X^{-1}A, \quad X = \omega \Lambda_{A}, 
\end{align*}
where $\Lambda_{A}$ denotes the matrix containing the diagonal elements of our matrix $A$. In general, for symmetric elliptic problems such as the Poisson equation, the optimal relaxation parameter is set at $\omega = \frac{3}{2}$. {\color{black}For the indefinite nonsymmetric Helmholtz problem, more tailored values of relaxation parameters should be used, especially outside the context of using multigrid as a preconditioner \cite{ernst2013multigrid, hadley2005complex}}.

To develop the theory, we start with the assumption that we are only using a post-smoother. In this case, the two-grid iteration matrix is given by 
\begin{align}
    T_0 = \left( I - P'{A_c}^{-1}PA \right) \left( I - X^{-1}A \right). \label{tgi1}
\end{align}

Generally the multigrid iteration scheme converges if the spectral radius of the two-grid iteration matrix is less than 1. Thus, the straightforward way to prove convergence is by showing that the spectral radius of the two-grid iteration matrix is strictly less than zero. If Dirichlet conditions are used, this naturally boils down to using Rigorous Fourier Analysis (RFA) to find analytical expressions for the eigenvalues of $T_0$. However, once we use Sommerfeld boundary conditions, this is no longer an option and we have to find another way to assess the boundedness of the spectral radius, {\color{black} which will be the focus of the next section, where we develop a general theoretical framework to assess the convergence, irrespective of the boundary conditions used.}

\subsection{Convergence} 
% We start by defining the two-grid operator in \cref{def1}. 
% \begin{def}[Two-grid operator] \label{def1} 
% Let $A \in \mathbb{C}^{n \times n}$ be a non-singular and non-defective matrix. $P,P' \in \mathbb{R}^{n \times m}$ are real valued prolongation and restriction operators. The $\omega-$Jacobi smoother is given by $X = \omega \Lambda_{A}$, where $\Lambda$ denotes the diagonal matrix containing the diagonal entries of $A$. If we allow a post-smoothing step, the two-grid operator is defined by 
% \begin{align}
%     T_0 = \left(I - P'{A_c}^{-1}PA \right) \left( I - X^{-1}A\right),
% \end{align}
% where the Galerkin condition is used to construct $A_{c} = P'AP \in \mathbb{C}^{m \times m}$.
% \end{def}
For convergence of the two-grid method we require that $\left \| T_{0} \right \|_{2} < 1$ independent of $h$. We can write $T_0$ as $T_0 = I - DA$, where $D$ represents an approximation to $A^{-1}$. This is shown in \cref{mgl1}. We use this to show that if the two-grid operator $T_0$ can be written in this form and the two-norm is bounded by one, then the multigrid method will converge.

\begin{lemma}{Convergence - I}\label{mgl1}
Let $A \in \mathbb{C}^{n \times n}$ be a non-singular and non-defective matrix, and let the two-grid operator be given by $T_0 = \left(I - P'{A_c}^{-1}PA \right) \left( I - X^{-1}A\right)$. Then, $T_{0}$ can be written as $T_{0} = I - DA$, with
\begin{align}
D = X^{-1} + P'{A_c}^{-1}P  - P'{A_c}^{-1}PX^{-1}    
\end{align}
Moreover we have ${T_0}^{H} = I - A^{H}D^{H}.$
\end{lemma}
\begin{proof}
Expanding $T_0$, and factoring in terms of $A$ we have
\begin{align}
    T_{0} &= I - X^{-1}A - P'{A_c}^{-1}PA  + P'{A_c}^{-1}PX^{-1}A, \\
    &= I - \left( X^{-1} + P'{A_c}^{-1}P  - P'{A_c}^{-1}PX^{-1} \right) A,\\
    &= I - DA, \label{p11}
\end{align}
where $D = X^{-1} + P'{A_c}^{-1}P  - P'{A_c}^{-1}PX^{-1}$. We similarly obtain
\begin{align}
    {T_0}^{H} &= {\left( \left(I - P'{A_c}^{-1}PA \right) \left( I - X^{-1}A\right) \right)}^{H}, \\
    &= {\left( I - X^{-1}A\right)}^{H}{\left( I - P'{A_c}^{-1}PA \right) }^{H}, \\
    &= \left( I - A^{H}X^{-1}\right) \left( I - A^{H}P'{\left({A_c}^{H}\right)}^{-1}P \right), \\
    &= I - A^{H}X^{-1} - A^{H}P'{\left({A_c}^{H}\right)}^{-1}P + A^{H}X^{-1}A^{H}P'{\left({A_c}^{H}\right)}^{-1}P, \\
    &= I - A^{H} \left( X^{-1} + P'{\left({A_c}^{H}\right)}^{-1}P - X^{-1}A^{H}P'{\left({A_c}^{H}\right)}^{-1}P \right) , \\
    &= I - A^{H}D^{H}, \label{p12}
\end{align}
where $D^{H} = X^{-1} + P'{\left({A_c}^{H}\right)}^{-1}P  - X^{-1}P'{\left({A_c}^{H}\right)}^{-1}P$.
\end{proof}

It has been shown that convergence of the two-grid operator implies multigrid convergence, see \cite{hackbusch2013multi}. The proof is based on an induction argument, where for the $W-$cycle, it is shown that the multigrid iteration matrix can be written as the two-grid iteration matrix plus a perturbation term. The perturbation term is bounded if the prolongation operator is a bounded linear operator and the smoother can be bounded by a constant. Note that a convergent smoothing scheme is not required. In fact, the smoothing property also holds for non-convergent iterations \cite{hackbusch2013multi}. It is only required that the error is smoothed up to a certain $\nu$, where $\nu$ is the number of smoothing steps. In other words, there exists a $\nu$ such that the error is reduced. It is not imperative that the smoothing property holds for $\nu = \infty$. 
We now work towards conditions for which $\left \| T_{0} \right \|_{2} < 1$ holds.

\begin{theorem}[Hermitian Positive Definiteness (HPD) implies Convergence] \label{mgt1} \\
%\begin{theorem}{Convergence - I}{mgt1} %\label{t0}
Let $A \in \mathbb{C}^{n \times n}$ be a non-singular and non-defective matrix. Let the two-grid operator be given  by $T_0 = \left(I - P'{A_c}^{-1}PA \right) \left( I - X^{-1}A\right)$, where $X$ is the $\omega-$Jacobi smoother, $P,P'\in \mathbb{R}^{n \times m}$ are the prolongation and restriction operator and $A_c = P'AP$. Let 
\begin{align}
    \Gamma = A^{H}D^{H} + DA - A^{H}D^{H}DA,
\end{align}
then ${T_0}^{H} {T_0} = I - \Gamma$. If $\Gamma$ is  positive definite, then $ {\lambda_j}\left( \Gamma \right) \in (0,2)$ and $\left \| T_{0} \right \|_{2} < 1$ independent of $h$. 
%\end{theorem}
\end{theorem}
\begin{proof}
Using equation \cref{mgl1}, equation \cref{p11} and \cref{p12}, we can write
\begin{align}
   {T_0}^{H} {T_0} &= \left( I - A^{H}D^{H} \right) \left( I - DA \right), \\
   &= I - \left( DA + A^{H}D^{H} - A^{H}D^{H}DA\right), \\
   &= I - \Gamma. \label{p16}
\end{align}
We can show a stronger result than just positive definiteness since $\Gamma$ is also Hermitian. To see this, note
\begin{align}
    {\Gamma}^{H} &= {\left( A^{H}D^{H} + DA - A^{H}D^{H}DA\right)}^{H}, \\
    &= {\left( DA + A^{H}D^{H} \right)}^{H} - {\left( A^{H}D^{H}DA\right)}^{H}, \\
    & = {\left( DA \right)}^{H} + {\left( A^{H}D^{H} \right)}^{H} - {\left( A^{H}D^{H}DA\right)}^{H}, \\
    &= A^{H}D^{H} + DA - A^{H}D^{H}DA = \Gamma. \label{p13}
\end{align}
Given that $\Gamma$  is HPD, the singular value decomposition (SVD) and eigendecomposition coincide, so we can find a unitary matrix $U\in \mathbb{C}^{n \times n}$ such that $\Gamma = U \Sigma_{\Gamma} U^{H}$, where $\Sigma$ is the diagonal matrix containing the singular resp. eigenvalues $(\sigma_j)$, where we have $\sigma_1 \geq \sigma_2 \geq \sigma_3 \hdots \geq \sigma_n > 0$. We have
\begin{align}
    \Gamma &= A^{H}D^{H} + DA - A^{H}D^{H}DA, %\\
    %& = U {\Sigma_{\Gamma}} U^{H} + \left( U {\Sigma_{\Gamma}} U^{H} \right)^{H} - \left( U {\Sigma_{\Gamma}} U^{H} \right)^{H}\left( U {\Sigma_{\Gamma}} U^{H} \right), %\\
%    & = U {\Sigma_{\Gamma}}  U^{H} + U{\Sigma_{\Gamma}}^{H}U^{H} - U{\Sigma_{\Gamma}} I{\Sigma_{\Gamma}}^{H}U^{H}, \label{p25}
\end{align}
which is equivalent to
\begin{align}
    \Sigma_{\Gamma} &= U A^{H}D^{H}U^{H} + UDAU^{H} - UA^{H}D^{H}DAU^{H}.\label{pp}
\end{align}
By our assumption on $\Gamma$ being positive definite, we have that for all $x \in \mathbb{C}^{n} \setminus\{ 0\}$,
\begin{align}
0 < x^{H}\Gamma x  &\Leftrightarrow 0 < x^{H}A^{H}D^{H}x + x^{H}DAx - x^{H}A^{H}D^{H}DAx. \label{eq143}%\\
% & \Leftrightarrow 0 < x^{H} U {\Sigma_{\Gamma}}  U^{H}x + x^{H} U{\Sigma_{\Gamma}}^{H}U^{H}x - x^{H}U{\Sigma_{\Gamma}} I{\Sigma_{\Gamma}}^{H}U^{H}x, \\
% & \Leftrightarrow 0 < x^{H}U{\Sigma_{\Gamma}} I{\Sigma_{\Gamma}}^{H}U^{H}x < x^{H} U {\Sigma_{\Gamma}}  U^{H}x + x^{H} U{\Sigma_{\Gamma}}^{H}U^{H}x.
\end{align}
Similarly, by using \cref{pp}, for all $x \in \mathbb{C}^{n} \setminus\{ 0\}$,
\begin{align}
0 < x^{H}\Sigma_{\Gamma} x  &\Leftrightarrow 0 < x^{H}U A^{H}D^{H}U^{H}x + x^{H}UDAU^{H}x - x^{H}UA^{H}D^{H}DAU^{H}x, \\
& \Leftrightarrow 0 < x^{H}UA^{H}D^{H}DAU^{H}x < x^{H}U A^{H}D^{H}U^{H}x + x^{H}UDAU^{H}x. \label{ppt}
\end{align}
Using the unitary invariance of the SVD, \cref{ppt} leads to 
\begin{align}
0 &< x^{H}UA^{H}D^{H}DAU^{H}x < x^{H}U A^{H}D^{H}U^{H}x + x^{H}UDAU^{H}x, \\
&\Rightarrow 0 < {\sigma_{j}(DA)}^2 < {\sigma_{j}(A^{H}D^{H} + DA)} \leq 2 {\sigma_{j}(DA)}, \\
& \Rightarrow \sigma_j(DA) < 2. \label{eq143-3}
\end{align}
This gives an upperbound for $\sigma_j(DA)$. We now develop the upper bound for the largest singular value of $\Gamma$, which is denoted by $\sigma_{1}(\Gamma)$. We know that $ 0 \leq \sigma_j(DA) < 2$. The case where $ 0 \leq \sigma_j(DA) < 1$ is trivial as that automatically renders $\left \| T_{0} \right \|_{2} < 1$. {\color{black}We focus on the case where $ 1 \leq \sigma_j(DA) < 2$. Recall, we assumed $\Gamma$ is HPD, but $DA$ does not necessarily have to be. As a result, the singular value decomposition for $\Gamma$ with unitary matrix $U$ might not coincide with the singular value decomposition for $DA$. 

Suppose, we can write $DA$ as $DA = W\Gamma_{DA}V^{H}$, where $WV^{H} \neq I$. We then still have
\begin{align}
    A^{H}D^{H}DA &= V {\Gamma_{DA}}^{H}W^{H}W\Gamma_{DA}V^{H}, \\
    V^{H}A^{H}D^{H}DAV &= {\Gamma_{DA}}^{H}\Gamma_{DA} = {\left(\Gamma_{DA}\right)}^2. 
\end{align}
Moreover, we also have for $v_j \in V$
\begin{align}
    DAv_j &= \sigma_j(DA)w_j \mbox{ and } \\
    \left( DAv_j \right)^{H} &= {v_j}^{H}A^{H}D^{H} = {w_j}^{H}\sigma_j(DA).
\end{align}
For $WV^{H} \neq I$ we will always have that for some index $j$ and $v_j \in V$ and $w_j \in W$
\begin{align}
    {v_j}^{H}w_j + {w_j}^{H}v_j < 2, 
\end{align}
where we used that the complex parts cancel out. We therefore have
\begin{align}
   v^{H}DAv + v^{H}A^{H}D^{H}v  =  {v_j}^{H}\sigma_j(DA)w_j +  {w_j}^{H}\sigma_j(DA)v_j < 2\sigma_j(DA). \label{eq143-2}
\end{align}

Taking $x = v_j$ in \cref{eq143} and using \cref{eq143-2}, $x^{H}\Gamma x$ can be bounded by
\begin{align}
 0 &< {v_j}^{H}\Gamma v_j = {v_j}^{H}A^{H}D^{H}v_j + {v_j}^{H}DAv_j - {v_j}^{H}A^{H}D^{H}DAv_j, \\ 
 &= \sigma_j(DA)\left( {v_j}^{H}w_j + {w_j}^{H}v_j  \right) - {\sigma_j(DA)}^2, \\
 &< 2\sigma_j(DA) - {\sigma_j(DA)}^2 = \sigma_j(DA)\left(2 - \sigma_j(DA)\right) , \\
 &< \sigma_j(DA) < 2.
\end{align}
}
%{\color{red} So we use the bounds we got from $\sigma_j(DA)$ to bound $\sigma_j(\Gamma)$, which coincides with $\lambda_j(\Gamma)$, given that this holds for all $j$.}
% % By rearranging \cref{pp} and evaluating in $j = 1$ (maximum), we have
% % \begin{align}
% %     \sigma_{1}( A^{H}D^{H} + DA  - \Sigma_{\Gamma}) &= \sigma_{1}(A^{H}D^{H}DAU^{H}) = {\sigma_{1}(DA)}^2. \label{pp1}
% % \end{align}
% %Note that we always have
% \begin{align}
%  \sigma_{1}(A^{H}D^{H} + DA)  - \sigma_{1}(\Gamma) \leq 2\sigma_{1}(DA)  - \sigma_{1}(\Gamma) 
% \end{align}
% % We can find two lower bounds for the left hand side of \cref{pp1}, which is given by
% % \begin{align}
% %  \sigma_{1}(A^{H}D^{H} + DA)  - \sigma_{1}(\Gamma) \leq 2\sigma_{1}(DA)  - \sigma_{1}(\Gamma) \leq  \sigma_{1}( A^{H}D^{H} + DA  - \Sigma_{\Gamma}).  
% % \end{align}
% Next, we claim that we must have
% \begin{align}
%  \sigma_{1}(A^{H}D^{H} + DA)  - \sigma_{1}(\Gamma) &\leq 2\sigma_{1}(DA)  - \sigma_{1}(\Gamma) \leq  {\sigma_{1}(DA)}^2. \label{ppx}
% \end{align}
% To see this, suppose the claim is false. Then for the right inequality we must have
% \begin{align}
%  2\sigma_{1}(DA)  - \sigma_{1}(\Gamma) > {\sigma_{1}(DA)}^2, \mbox{ with $0 \leq \sigma_j(DA) < 2$}. \label{pp2}
% \end{align}
% \cref{pp2} can only be true 

% for this range of $\sigma_{1}(DA)$ in terms of $\sigma_{1}(\Gamma)$ gives the condition $\sigma_{1}(\Gamma) \leq 0$, which contradicts our assumption of $\Gamma$ being positive definite. In fact for \cref{ppx} to be valid for the entire range of $1 \leq \sigma_{1}(DA) < 2$, we must have 
% $1 \leq \sigma_{1}(\Gamma) < 2$. 

As $\Gamma$ is positive definite, the eigenvalues and singular values coincide and we have 
\begin{align}
    \lambda_{\min}(\Gamma) > 0 \mbox{ and } \lambda_{\max}(\Gamma) < 2. \label{pfin}
\end{align}
We are now ready to bound the two-grid operator $T_0$ using the 2-norm over $\mathbb{C}^{n \times n}$
\begin{align}
    \| {T_0} \|_{2} &= \sup_{x \in \mathbb{C}^{n},  \|x\|_{2} = 1} { \| \left( {T_0}^{H} {T_0 } \right) x \|_{2} } = \sqrt{\rho \left( {T_0}^{H} {T_0}  \right )}, \nonumber \\
    &= \sqrt{\rho \left( I - \left( A^{H}D^{H} + DA - A^{H}D^{H}DA \right) \right )}, \\
    &= \sqrt{\max_{1\leq j \leq n} | \lambda_{j} \left( I - \left( A^{H}D^{H} + DA - A^{H}D^{H}DA \right) \right) | }, \\
    & \leq \sqrt{| 1 - \min_{1\leq j \leq n} \lambda_{j} \left( A^{H}D^{H} + DA - A^{H}D^{H}DA \right) | } , \label{p23x} \\
    & < 1 \nonumber,
\end{align}
where we used that trivially $0 < \min \lambda_{j} \left( A^{H}D^{H} + DA - A^{H}D^{H}DA \right)$ by positive definiteness of $\Gamma$. As regards, the largest eigenvalue $\max \lambda_{j} \left( \Gamma \right) < 2$, we obtain
\begin{align}
    \| {T_0} \|_{2} &= \sup_{x \in \mathbb{C}^{n},  \|x\|_{2} = 1} { \| \left( {T_0}^{H} {T_0 } \right) x \|_{2} } = \sqrt{\rho \left( {T_0}^{H} {T_0}  \right )}, \nonumber \\
    &= \sqrt{\rho \left( I - \left( A^{H}D^{H} + DA - A^{H}D^{H}DA \right) \right )}, \\
    &= \sqrt{\max_{1\leq j \leq n} | \lambda_{j} \left( I - \left( A^{H}D^{H} + DA - A^{H}D^{H}DA \right) \right) | }, \\
    & \leq \sqrt{| 1 - \max_{1\leq j \leq n} \lambda_{j} \left( A^{H}D^{H} + DA - A^{H}D^{H}DA \right) | } , \label{p23y} \\
    & < \sqrt{|1 - 2|} = 1, \mbox{ using \cref{pfin}.}
\end{align}
Thus, if $\Gamma$ is positive definite then $\lambda_j(\Gamma) \in (0,2)$ and consequently $\lambda_j(T_0) \in (0,1)$. 
\end{proof} 

{\color{black}\cref{mgt1} shows that if $\Gamma$ is positive definite, then the eigenvalues of $T_0$ will lie in the interval $(0,1)$ which automatically leads to a bound for the 2-norm of the two-grid operator. Note that this is independent of $h$, as no constant appears in \cref{pfin} which depends on $h$. Now that we know that a sufficient conditions for convergence is to have $\Gamma$ be positive definite, we need a way to practically assess this. One approach is by simplifying the expression for $\Gamma$ and finding another sufficient condition for positive definiteness, which is more easy to analyze. This will be constructed in \cref{mgc1}}. We again start by defining the simplified operator in \cref{def2s}. 
\begin{definition}[Simplified operator] \label{def2s} 
Let $A \in \mathbb{C}^{n \times n}$ be a non-singular and non-defective matrix. $P,P' \in \mathbb{R}^{n \times m}$ are real valued prolongation and restriction operators and $A_{c} = P'AP$. The $\omega-$Jacobi smoother is given by $X = \omega \Lambda_{A}$, where $\Lambda$ denotes the diagonal matrix containing the diagonal entries of $A$. Then we define $\tilde{D} = X^{-1} + P'{A_c}^{-1}P$ such that $\tilde{\Gamma} = A^{H}\tilde{D}^{H} + \tilde{D}A - A^{H}\tilde{D}^{H}\tilde{D}A$. 
\end{definition}
Note that in the definition given above, $\tilde{D}$ does not contain the term $P'{A_c}^{-1}PX^{-1}A$, which is part of the original operator $D$. %as given by def1 of T0

\begin{corollary}{Convergence - II}\label{mgc1}
Let $A \in \mathbb{C}^{n \times n}$ be a non-singular and non-defective matrix. Let $\tilde{\Gamma}$ be defined as in \cref{def2s}. If $\tilde{\Gamma}$ is positive definite, then $\Gamma$ from \cref{mgt1} is positive definite. 
 \end{corollary}
\begin{proof}
To see this, we write $\Gamma$ as
\begin{align}
    \Gamma &= A^{H}D^{H} + DA - A^{H}D^{H}DA, \\
    &= \left(A^{H}\tilde{D}^{H} + \tilde{D}A - A^{H}\tilde{D}^{H}\tilde{D}A \right) + {A^H}{X^{-1}P'{{\left({A_c}^{H}\right)}^{-1}}P}{P'{A_c}^{-1}PX^{-1}}A\\
    &- A^{H}\tilde{D}^{H}P'{A_c}^{-1}PX^{-1}A - {A^H}X^{-1}P'{{\left({A_c}^{H}\right)}^{-1}}P\tilde{D}A.
\end{align}
We assume $\tilde{\Gamma}$ is positive definite and HPD respectively. Now, assuming the opposite, i.e. $\Gamma$ is not HPD implies that $\exists x \in \mathbb{C}^{n} \setminus\ 	\{ 0 \}$ such that
\begin{align*}
&x^{H}\tilde{\Gamma}x + {x^{H}}{A^H}{X^{-1}P'{{\left({A_c}^{H}\right)}^{-1}}P}{P'{A_c}^{-1}PX^{-1}}Ax \hspace{2mm} (\ast)\\
&- {x^{H}}A^{H}\tilde{D}^{H}P'{A_c}^{-1}PX^{-1}Ax - {x^{H}}{A^H}X^{-1}P'{{\left({A_c}^{H}\right)}^{-1}}P\tilde{D}Ax \leq 0. \hspace{2mm} (\ast \ast)
\end{align*}
By assumption and the fact that the second term in $(\ast)$ contains a quadratic form, $(\ast)$ is always positive for any $x \in \mathbb{C}^{n} \setminus\ 	\{ 0 \}$. But if $\Gamma$ is not positive definite then using $(\ast \ast)$ we must also have
\begin{align}
    0 < x^{H}\tilde{\Gamma}x + {x^{H}}{A^H}{X^{-1}P'{{\left({A_c}^{H}\right)}^{-1}}P}{P'{A_c}^{-1}PX^{-1}}Ax, \\
    \leq {x^{H}}A^{H}\tilde{D}^{H}P'{A_c}^{-1}PX^{-1}Ax + {x^{H}}{A^H}X^{-1}P'{{\left({A_c}^{H}\right)}^{-1}}P\tilde{D}Ax \label{newtei1}
\end{align}

{\color{black}
We show that $P'{A_c}^{-1}PX^{-1}A$ has singular values equal to zero. Observe that the rank of $P'{A_c}^{-1}PX^{-1}A$ is $\frac{n}{2}$. Suppose we can write the singular value decomposition as $\Hat{\Gamma} = W^{H}\left(P'{A_c}^{-1}PX^{-1}A\right)V$. In this case we have $V$ and $W$ to represent two unitary complex matrices as we have not assumed that $P'{A_c}^{-1}PX^{-1}A$ is HPD. We know that $P'{A_c}^{-1}P$ is singular and half of all singular values are equal to zero. We thus need to show that $P'{A_c}^{-1}PX^{-1}A$ has zero singular values as well. In order to do this, we use an analogous version of Ostrowski’s theorem. Here instead of using eigenvalues, we use singular values. According to the theorem, for $K \in \mathbb{C}^{n \times n}$ and $Y \in \mathbb{C}^{n \times n}$ with $Y$ is non-singular, we have 
\begin{align}
    \sigma_j(KY) \leq \sigma_1(Y)\sigma_j(A), \quad j = 1,2, \hdots n,  \label{newtei}
\end{align}
where $\sigma_1 \geq \sigma_2 \geq \hdots \sigma_n$. Applying \cref{newtei} to our case by taking $K = P'{A_c}^{-1}P$ and $Y = X^{-1}A$, we obtain
\begin{align}
    \sigma_j\left( P'{A_c}^{-1}P X^{-1}A \right) \leq \sigma_1(X^{-1}A)\sigma_j(P'{A_c}^{-1}P), \quad j = 1,2, \hdots n.  \label{newtei2}
\end{align}
We know that $\sigma_1(X^{-1}A) > 0$ given that both $X$ and $A$ are non-singular. However, we know that $P'{A_c}^{-1}P$ is singular and has zero singular values. Thus there exists indices $\tilde{j}$ such that 
\begin{align}
    \sigma_{\tilde{j}}\left( P'{A_c}^{-1}P X^{-1}A \right) \leq 0, \quad \tilde{j} \in \{ 1,2, \hdots n \} . \label{newtei3}
\end{align}
Trivially  $0 \leq \sigma_{\tilde{j}}\left( P'{A_c}^{-1}P X^{-1}A \right)$ is a lower bound for the smallest singular value. Combining this with \cref{newtei3} leads to
\begin{align}
   0 \leq \sigma_{\tilde{j}}\left( P'{A_c}^{-1}P X^{-1}A \right) \leq 0, \quad 	\exists \tilde{j} \in \{ 1,2, \hdots n \}, \label{newtei4}
\end{align}
and so by the interlacing theorem, we must have that there exist at least one singular value at an index $\tilde{j}$ such that $\sigma_{\tilde{j}}\left( P'{A_c}^{-1}P X^{-1}A \right) = 0$. Naturally, the eigenvalue of $P'{A_c}^{-1}P X^{-1}A$ corresponding to the index $\tilde{j}$ must be zero as well, given that the singular value is zero. 

We know that for any complex matrix, there exists a unitary matrix such that we can write it in the Schur decomposition. Without the assumption of normality, it solves a general eigenvalue problem of the form $Kv = \lambda Bv$ for some matrix $B$. 

We take $x = v_{\tilde{j}}$, corresponding to the zero eigenvalue, which we denote by $\lambda_{\tilde{j}}^{0}$ such that
\begin{align}
P'{A_c}^{-1}P X^{-1}Av_{\tilde{j}} = \lambda_{\tilde{j}}(P'{A_c}^{-1}P X^{-1}A)Bv_{\tilde{j}} = 0Bv_{\tilde{j}}. \label{newtei5}
\end{align}
Substituting \cref{newtei5} into \cref{newtei1} gives
\begin{align}
    0 &< {v_{\tilde{j}}}^{H}\tilde{\Gamma}v_{\tilde{j}} + {{v_{\tilde{j}}}^{H}}{A^H}{X^{-1}P'{{\left({A_c}^{H}\right)}^{-1}}P}{P'{A_c}^{-1}PX^{-1}}Av_{\tilde{j}}, \label{newtei60} \\
    &\leq {{v_{\tilde{j}}}^{H}}A^{H}\tilde{D}^{H}P'{A_c}^{-1}PX^{-1}Av_{\tilde{j}} + {{v_{\tilde{j}}}^{H}}{A^H}X^{-1}P'{{\left({A_c}^{H}\right)}^{-1}}P\tilde{D}Av_{\tilde{j}}, \label{newtei6} \\
    &= \lambda_{\tilde{j}}^{0}{{v_{\tilde{j}}}^{H}}A^{H}\tilde{D}^{H}B{v_{\tilde{j}}} + \lambda_{\tilde{j}}^{0} {{v_{\tilde{j}}}^{H}} B^{H}\tilde{D}A {v_{\tilde{j}}} = 0.  
\end{align}
} 
This gives a contradiction as $\tilde{\Gamma}$ is assumed to be positive definite and must always be larger than zero for any $x \in \mathbb{C}^n \setminus\ 	\{ 0 \}$. We thus have that if $\tilde{\Gamma}$ is positive definite, then $x^{H}\Gamma x > 0 \mbox{ }\forall x \in \mathbb{C}^{n} \setminus\ 	\{ 0 \}$, i.e. $\Gamma$ must be positive definite as well. 
\end{proof} 

% We take $x = w_j$, where $w_j$ is an right normalized eigenvector of $P'{A_c}^{-1}PX^{-1}A$ corresponding to a zero-eigenvalue, which we denote by ${\lambda_j}^{0} = 0$. It is easy to see that, while $A_c$ is non-singular, $P'{A_c}^{-1}PX^{-1}A$ is not. \\
% For example, we can bound the eigenvalues of $P'{A_c}^{-1}PX^{-1}A$ from below by\\ $\min{X^{-1}_{jj}}IP'{A_c}^{-1}PA$ and from above by $\max{X^{-1}_{jj}}IP'{A_c}^{-1}PA$. $\min{X^{-1}_{jj}}$ and $\max{X^{-1}_{jj}}$ are now constants, and we know that $P'{A_c}^{-1}PA$ has exactly half of its eigenvalues equal to one, and the other half equal to zero. We therefore know that $P'{A_c}^{-1}PX^{-1}A$ must have zero eigenvalues as well. Consequently, we obtain 
% \begin{align}
%     0 &< {w_j}^{H}\tilde{\Gamma}{w_j} + {{w_j}^{H}}{A^H}{{X}^{-1}P'{{\left({A_c}^{H}\right)}^{-1}}P}{P'{A_c}^{-1}P{X}^{-1}}A{w_j}, \\
%     &\leq {{w_j}^{H}}A^{H}\tilde{D}^{H}P'{A_c}^{-1}P{X}^{-1}A{w_j} + {{w_j}^{H}}{A^H}{X}^{-1}P'{{\left({A_c}^{H}\right)}^{-1}}P\tilde{D}A{w_j}, \\
%     &= {\lambda_j}^{0}{{w_j}^{H}}A^{H}\tilde{D}^{H}{w_j} + {\bar{\lambda_j}^{0}}{{w_j}^{H}}\tilde{D}A{w_j} = 0.
% \end{align}

Using the above, we now have that it is sufficient to check if $\tilde{\Gamma}$ is positive definite. We therefore test for convergence by verifying whether $\tilde{\Gamma}$ is positive definite. We use the Cholesky decomposition to determine whether $\tilde{\Gamma}$ is positive definite, as any HPD system matrix $B$ can be written as $B = LL^{\ast}$, where $L \in \mathbb{C}^{n \times n}$ is a lower-triangular matrix. Another method to check for positive definiteness without having to compute the Cholesky factors is by checking all of the following
\begin{enumerate}
    \item $b_{ii} > 0 \forall i$
    \item $b_{ii} + b_{jj} > 2 \left|\Re \left[ b_{ij} \right] \right| \mbox{ for } i \neq j$
    \item The element with the largest modulus lies on the diagonal
    \item $\det (B) > 0$
\end{enumerate}

{\color{black}The results from this section show that if we have an HPD ${\Gamma}$ matrix, then the two-grid method converges. This can also be tested for a simplified operator $\tilde{\Gamma}$. Note that so far, all theoretical considerations are for general matrices. Next, we will apply this theory to the Helmholtz operator.}

\section{Convergence for indefinite nonsymmetric Helmholtz operators} \label{sec3} %% 12-5 \newline
{\color{black}The previous section has stipulated properties which we can check in order to find a convergent configuration for the indefinite nonsymmetric Helmholtz equation. For this we need to consider optimality conditions for the smoother and coarse grid correction scheme.

Starting with the smoother, we will the use classical $\omega-$Jacobi smoother due to its natural parallelizability. Even though the smoother diverges for Helmholtz problems, we will show that combined with CSL-coarsening and higher-order prolongation and restriction, the $\tilde{\Gamma}$ and consequently $\Gamma$ operator becomes HPD, leading to a convergent solver on a complete V-cycle hiearchy. In the remainder of this paper, unless stated otherwise, we will use aggregate values from the literature presented for the two-step Jacobi as our $\omega$ in the weighted-Jacobi method, which boils down to keeping $\omega = 4.5$. Our own experimental results in \cref{opom} are in line with the aggregate values from the literature, confirming relaxation parameters which are generally lower than the underrelaxation parameters reported for symmetric elliptic problems, such as Poisson's equation.

}
%{\color{purple}
% In this section, we prove two properties which we need in order to construct a robust solver: convergence and an optimal smoother. We start with convergence and then work our way towards deriving optimality conditions for the $\omega-$Jacobi smoother. The $\omega-$Jacobi smoother is known to diverge for Helmholtz problems. However, we show that with the right conditions, we can obtain a convergent solver. 

% [{\color{purple} hier iets schrijven over de relaxation parameter. Generally underrelaxed, omega tussen 0 en 1, voor poisson is dat 2/3, verwijzingen linz, sheikh, erlangga, binzubair. Hier alvast introduceren dat het voor helmholtz andere waarde kan aannemen en dan misschien verwijzen naar optimality voor Helmholtz hoofdstuk. Die inleding is logischer. Tabel voor Helmholtz ook eventueel verplaatsen dan.}] }

{\color{black}We proceed by evaluating these properties for the Helmholtz operator. We use MP2-A with constant wavenumber $k$. We numerically investigate the changes we can make to the design of the coarse-grid system in order to achieve convergence. In this work, we consider two options; we either use higher-order interpolation schemes to construct the prolongation and restriction operator and/or we use the CSL, which will be denoted by $C$, as the system on which we apply the coarsening operations. We use the higher-order interpolation scheme based on the quadratic rational Bezi\'er curve from \cite{dwarka2020scalable,dwarka2022scalable}. There, the scheme was used to construct the deflation vectors for a two-level and multi-level deflation preconditioner, which showed wavenumber independent convergence for the two-level scheme and close to wavenumber independent convergence for the multi-level scheme. The main motivation for studying these type of higher-order schemes is that they have been shown to reduce the projection error \cite{dwarka2020scalable}. 
Using this scheme, the prolongation operator acts on a grid function as follows}
\begin{equation}
P\left[u_{2h}  \right ]_{i} = \left\{\begin{matrix}
\begin{matrix}
\frac{1}{8}  \left(\left[u_{2h}  \right ]_{\left(i - 2\right) / 2} + 6 \left[u_{2h}  \right ]_{\left(i \right) / 2} + \left[u_{2h}  \right ]_{\left(i + 2\right) / 2} \right) & \mbox{if $i$ is even,}\\
\frac{1}{2}  \left(\left[u_{2h}  \right ]_{\left(i - 1\right) / 2} + \left[u_{2h}  \right ]_{\left(i + 1\right) / 2} \right) &  \mbox{if $i$ is odd,}
\end{matrix}
\end{matrix}\right\}, \label{xex}
\end{equation}
for $ i = 1, \dots, N$ and for $i = 1, \dots, \frac{N}{2}$, where $N$ denotes the size of the fine-level linear system. Note that the second line ($i$ = odd) is in fact the linear interpolation scheme. Thus, we now have a combination of both a higher-order and linear interpolation scheme for the nodes $i$ within the numerical domain. 

The inclusion of the complex shift was first applied to multigrid schemes in \cite{cools2014new}. There, the main motivation relies on spectral analysis of the two-grid operator and the observation that the coarse-grid eigenvalues can never approach zero due to the complex part. We proceed by a different yet general theoretical argument, which is independent of the boundary conditions. Moreover, we include a constant complex shift compared to a level-dependent one and also evaluate the effect of a classic smoother compared to the polynomial GMRES(3) smoother.  

In \cref{conv1}, we show that the use of the CSL in combination with a higher-order interpolation scheme in fact leads to an HPD system. Note that while our conditions do not require the actual computation of the spectral norm, we will do so for sake of illustration and completeness. 
\begin{table}[!hbt]
\renewcommand{\arraystretch}{1.5}
\centering
\scalebox{1}
{
\begin{tabular}{ccc|cc}
\hline
& \multicolumn{2}{c}{Linear} & \multicolumn{2}{c}{Bezi\'er} \\ \hline
\multicolumn{1}{c|}{$k$} & $A$                 & $C$                 & $A$                  & $C$                 \\ \hline
\multicolumn{1}{c|}{$5$}  & $\times {\color{white}-} 2.284$     & $\times {\color{white}-} 1.304$    & $\checkmark {\color{white}-} 0.991$  & $\checkmark {\color{white}-} 0.911$               \\
\multicolumn{1}{c|}{$10$} & $\times {\color{white}-} 5.888$     & $\times {\color{white}-} 1.351$    & $\times {\color{white}-} 1.105$  & $\checkmark {\color{white}-} 0.913$               \\
\multicolumn{1}{c|}{$20$} & $\times {\color{white}-} 8.786$     & $\times {\color{white}-} 1.328$    & $\times {\color{white}-} 1.306$  & $\checkmark  {\color{white}-} 0.951$               \\
\multicolumn{1}{c|}{$30$} & $\times {\color{white}-} 10.660$    & $\times  {\color{white}-} 1.325$    & $\times {\color{white}-} 1.504$  & $\checkmark {\color{white}-} 0.984$ \\ \hline
\end{tabular}
}
\hspace{0.5cm}
\captionsetup{width=1 \textwidth}
\caption{At the right of each entry, the spectral norm of the two-grid operator is given. Linear uses linear interpolation to construct $P,P'$. Bezi\'er uses rational quadratic Bezi\'er interpolation. $A$ represents $A_c = P'AP$. $C$ represents $A_c = P'CP$, where $C$ denotes the CSL with complex shift $\beta_2 = 0.7$. In all cases, one post-smoothing step is used. Left of each entry, $\checkmark$ denotes that $\tilde{\Gamma}$ is HPD and  $\times$ denotes it is not.}
\label{conv1}
\end{table}
\FloatBarrier
Several interesting observations can be made. First of all, in almost call cases when using the classical smoothing operator, a combination of \textit{both} higher-order interpolation and CSL for coarsening is needed to obtain an HPD matrix and consequently a spectral norm which is bounded by one. Thus, using either level independent CSL-coarsening or higher-order prolongatrion and restriction, will not solve the problem.

\FloatBarrier

\subsection{Optimality}
{\color{black}
Apart from having a condition to check for convergence, we would like to have a measure of how fast we can expect the two-grid method to converge. By having a positive definite $\Gamma$, we know that the spectral norm is bounded from above by one. To obtain a more qualitative estimate, we study the condition number of $\Gamma$ and $\tilde{\Gamma}$. Using the condition numbers we obtain a sharper bound in \cref{mgt1}. This entails that the smallest eigenvalue of $\tilde{\Gamma}$ will be a lower bound to the smallest eigenvalue of $\Gamma$.}
%In particular if ${\| \tilde{\Gamma} \|} \leq {\| {{\Gamma}} \|}$ in any norm then ${\| \tilde{\Gamma}^{-1} \|} \leq {\| {{\Gamma}}^{-1} \|}$. }
\begin{corollary}{Convergence Optimality Condition - I}\label{mgc2}
%\begin{cor}[Convergence - III] \label{p21}
Let $A \in \mathbb{C}^{n \times n}$ be a non-singular and non-defective matrix and let $\tilde{D} = X^{-1} + P'{A_c}^{-1}P$. If $\tilde{\Gamma} = A^{H}\tilde{D}^{H} + \tilde{D}A - A^{H}D^{H}DA$ is positive definite and ${\| {{\Gamma}} \|} \leq  {\| \tilde{\Gamma} \|} $ together with $\kappa(\tilde{\Gamma}) \leq  \kappa(\Gamma)$ then $\lambda_{\min}(\tilde{\Gamma})$ is given and bounded by
\begin{align}
     \frac{{\| {\tilde{\Gamma}} \|}}{\kappa(\tilde{\Gamma})}  \leq \lambda_{\min}({\Gamma}). \label{p24}
\end{align}
 \end{corollary}
\begin{proof}
Note that $\tilde{\Gamma}$ and $\Gamma$ are HPD and so all eigenvalues are real, so we can consider ${\| {{\Gamma}}^{-1} \|}$ and ${\| {\tilde{\Gamma}}^{-1} \|}$. Moreover, $\forall \mbox{ matrices } B \in \mathbb{C}^{n \times n}, \rho(B) = \max_{1\leq j \leq n}|\lambda_j(B)| \leq {\| B \|}$ in any associative norm. If $\kappa(\tilde{\Gamma}) \leq  \kappa(\Gamma)$, then we obtain
\begin{align}
    {\| {{\Gamma}} \|} \leq {\| \tilde{\Gamma} \|} \Rightarrow \frac{1}{{\| {\tilde{\Gamma}} \|}} \leq  \frac{1}{{\| {\Gamma} \|}}  \Rightarrow \frac{\kappa(\tilde{\Gamma})}{{\| {\tilde{\Gamma}} \|}} \leq \frac{\kappa(\Gamma)}{{\| {\Gamma} \|}}   \Rightarrow \| \tilde{\Gamma}^{-1} \| \leq \| {\Gamma}^{-1} \|. \label{pxf}
\end{align}
Using that the systems are HPD, we have 
\begin{align}
\rho(\Gamma^{-1}) = \max_{1\leq j \leq n}|\lambda_j(\Gamma^{-1})| =  \min_{1\leq j \leq n}{|\lambda_j(\Gamma)|}^{-1} \leq \| {\Gamma}^{-1} \|.    
\end{align}
Note that this gives the reciprocal of the smallest eigenvalue of $\Gamma$. To get the actual smallest eigenvalue of $\Gamma$, we have to take the reciprocal, which provides us with inequality \cref{p24} and a lower bound to $\lambda_{\min}({\Gamma})$.
\end{proof}  

We can use \cref{mgc2} to construct a sharper bound in \cref{mgt1}. 
%This is stated in Corollary \ref{co:mgc3}.
%\begin{cor}[Convergence - IV] \label{p22}
\begin{corollary}{Convergence Optimality Condition - II}\label{mgc3}
Let $A \in \mathbb{C}^{n \times n}$ be a non-singular and non-defective indefinite matrix. Let $T_{0}$ be such that we can write $T_0$ as $T_{0} = I - DA$ with $D = X^{-1} + P'{A_c}^{-1}P  - P'{A_c}^{-1}PX^{-1}$. Let $\Gamma = A^{H}{D}^{H} + DA - A^{H}{D}^{H}{D}A$. If $\tilde{\Gamma} = A^{H}\tilde{D}^{H} + \tilde{D}A - A^{H}\tilde{D}^{H}\tilde{D}A$ with $\tilde{D} = X^{-1} + P'{A_c}^{-1}P$, is positive definite and ${\| \tilde{\Gamma} \|} \leq {\| {{\Gamma}} \|}$ together with $\kappa(\Gamma) \leq \kappa(\tilde{\Gamma})$, then
\begin{align*}
\|{T_0} \|_{2} < \sqrt{\left| 1 - \frac{\| \tilde{\Gamma} \|_{\text{}}}{\kappa_{\text{}}(\tilde{\Gamma})}\right| } < 1. 
\end{align*}
 \end{corollary}
\begin{proof}
Using \cref{mgt1} we can write ${T_{0}}^{H}T_0 = I - \Gamma$. Using \cref{p23x} from \cref{mgt1} and substituting \cref{p24} from \cref{mgc3} have
\begin{align}
    \| {T_0} \|_{2} &= \sup_{x \in \mathbb{C}^{n},  \|x\|_{2} = 1} { \| \left( {T_0}^{H} {T_0 } \right) x \|_{2} } = \sqrt{\rho \left( {T_0}^{H} {T_0}  \right )}, \nonumber \\
    & = \sqrt{| 1 - \min \lambda_{j} \left( \Gamma \right) | }, \\
    & < \sqrt{| 1 - \min \lambda_{j} \left( \tilde{\Gamma} \right) | }, \\
    & < \sqrt{\left| 1 - \frac{\| \tilde{\Gamma} \|_{\text{}}}{\kappa_{\text{}}(\tilde{\Gamma})}\right| } < 1. \label{msms}
    \end{align} \qed
\end{proof}

\subsection{Optimal convergence for Helmholtz} \label{opom}
From \cref{sec3}, we have observed that we need a clever combination of coarsening using the CSL, higher-order interpolation schemes and the right amount of smoothing steps in order to obtain HPD systems. {\color{black}The results from \cref{conv1}, however, indicate that convergence, in terms of total number of iterations, will be slow as the spectral norm is still close to one.} 
% We can use \cref{mgc3} to obtain some insights into how we can optimize the convergence by optimizing the remaining parameter we can work with; the relaxation parameter $\omega$ in the $\omega-$Jacobi smoother. We let $\omega$ vary and report the value of \cref{p24} from \cref{mgc2} in \cref{opt1}. Note that smaller these number, the closer to one the two-grid spectral norm will be, given \cref{msms}. {\color{red} Consequently, we want this value to be larger as the number of smoothing steps increases. If not, performing additional smoothing steps leads to deteriorating convergence for $\omega-$Jacobi, which is generally known to be divergent for indefinite problems. However, we reiterate that despite the smoother being divergent, it can stil exhibit smoothing properties and be part of a converging solver.} 
% {\color{blue} waarom we uiteindelijk kiezen voor 4.5, merk op dat van 1 smoothing stap naar 2, de bovengrens kleiner wordt voor 1.5 en dus idd naar mate je meer gaat smoothen gaat het closer naar 1. bij die overige blijft het 'stijgen' ipv dalen want hoe kleiner dit getal hoe slechter je wilt dus dat dit getal groeit met de smoothing steps $\nu$. Sws weten we dat eig(i - omega jacobi a) altijd groter dan 1 zullen zijn dus daarom kan je niet triviaal bounden met product van evals van coarse grid corrector en smoother matrix.}

\begin{table}[!hbt]
\centering
\scalebox{0.95}{
\begin{tabular}{c|cc|cc|cc|cc|cc}
& \multicolumn{2}{c|}{$\omega=1.5$} &  \multicolumn{2}{c|}{$\omega=2$} & \multicolumn{2}{c|}{$\omega=2.5$} &  \multicolumn{2}{c|}{$\omega=4.5$}  &  \multicolumn{2}{c}{$\omega=7$} \\ \hline
% & \multicolumn{2}{c|}{$N = 6724$} &  \multicolumn{2}{c|}{$N = 26244$} & \multicolumn{2}{c|}{$N = 57600$} &  \multicolumn{2}{c|}{$N = 102400$} &  \multicolumn{2}{c}{$N = 160000$} \\ 
% & \multicolumn{2}{c|}{$N_{D} = 8$} &  \multicolumn{2}{c|}{$N_{D} = 8$} & \multicolumn{2}{c|}{$N_{D} = 4$} &  \multicolumn{2}{c|}{$N_{D} =8$} &  \multicolumn{2}{c}{$N_{D} = 4$} \\ \hline
$k$ & $\nu = 1$   & $\nu = 2$ & $\nu = 1$   & $\nu = 2$   & $\nu = 1$    & $\nu = 2$   & $\nu = 1$    & $\nu = 2$   & $\nu = 1$   & $\nu = 2$    \\ \hline 
$5$  & $0.373$  &$0.413$ & $0.273$ & $0.405$ & $0.206$  & $0.389$ & $0.088$  & $0.200$ & $0.031$ & $0.123$   \\
$10$ & $0.137$  &$0.140$ & $0.128$ & $0.139$ & $0.112$  & $0.137$ & $0.065$  & $0.116$ & $0.028$ & $0.081$   \\
$20$ & $0.030$  &$0.028$ & $0.029$  & $0.029$ & $0.028$  & $0.030$ & $0.022$  & $0.028$ & $0.012$ & $0.025$   \\
$30$ & $0.011$  &$0.009$ & $0.011$  & $0.010$ & $0.010$  & $0.011$ & $0.008$  & $0.011$ & $0.001$ & $0.009$   \\
\end{tabular}}
\caption{We report the value of $\| \tilde{\Gamma} \|\kappa(\tilde{\Gamma})^{-1}$ in the $p = 1$ norm, with $\tilde{\Gamma}$ from \cref{mgc2}. $\nu$ denotes the number of {$\omega$-Jacobi} smoothing steps.}
%}
\label{opt1}
\end{table}

The results from \cref{opt1} reveal that a lower $\omega$ is in favor when we use one post-smoothing steps and the wavenumber is small. For $\omega = 1.5$, we observe that increasing $k$ and $\nu$ leads to an decrease in $\| \tilde{\Gamma} \| \kappa(\tilde{\Gamma}) ^{-1}$. This means that as we perform more smoothing steps, the two-grid spectral norm moves closer and closer to one. For convergence, we require the opposite: the addition of a few more smoothing steps should lead to an increase rather than a decrease in $\| \tilde{\Gamma} \| \kappa(\tilde{\Gamma}) ^{-1}$. 
Based on these results, the optimal estimate for $\omega$ lies between $[2,4.5]$ for the reported values of $k$. {\color{black} These values are different from the general consensus on the optimal relaxation parameter for poisson-type of problems. While we focus on the classical smoother in this work, our reported values do seem to coincide with underrelaxation parameters observed in other studies where a two-step weighted Jacobi smoother is analyzed \cite{hadley2005complex, ernst2013multigrid}}
% While in general the application of the $\omega-$Jacobi smoother diverges for highly indefinite Helmholtz problems, choosing the right $\omega$ combined with a higher-order interpolation scheme and alternative coarsening strategies can lead to a simple yet convergent two-grid scheme. 

% \subsection{Convergence of the multigrid cycle}
% From here
% \subsection{Optimality}
% {\color{red} Hier laten we zien dat de $\omega$ van $\omega-$Jacobi invloed heeft op de conditie $\kappa(\Gamma) \leq \kappa(\tilde{\Gamma})$. We willen dus een $\omega$ kiezen zodanig dat aan de conditie wordt voldaan maar het liefst is $\kappa(\tilde{\Gamma})$ zo klein mogelijk wordt. Dit gebeurt als we dus coarsening toepassen op CSLP, hogere orde prolongatie en restriction gebruiken en $\omega > 1.5$ kiezen. (max 1 pagina)}

\section{Numerical Experiments} \label{sec5} 
\FloatBarrier
In this section we examine the convergence behavior of the multigrid solver. We use the model problems from \cref{sec2}. Unless stated otherwise, we use 10 grid points per wavelength, which is equivalent to $kh = 0.625$. All experiments are implemented sequentially on a Dell laptop using 8GB RAM and a i7-8665U processor. In all experiments, we set the relative tolerance to $10^{-5}$. We allow for coarsening until the dimension of the underlying linear system $N_c$ is smaller than 10. The maximum size of the linear system on the coarsest grid is therefore $10 \times 10$. We use the $\omega-$Jacobi smoother with the approximate optimal relaxation parameter from \cref{opom}, i.e. $\omega = 4.5$ as we want to test for large $k$. Moreover, we perform the coarse-grid corrections using the higher-order prolongation and restriction operator, together with the CSL. For the latter, the complex shift is set at $\beta_2 = 0.7$ unless stated otherwise. 
% The scope and aim of these experiments are to assess whether we can reach convergence using the new set-up. We do not yet require the convergence to be wavenumber independent, as that will require the study of different smoothers. We aim to establish a baseline for the solver from which we can 

\subsection{2D Constant \texorpdfstring{$k$}{k}}
For the constant wavenumber problem (MP 2-A), we first start by confirming that we have $h-$independent convergence. Recall that the constant from \cref{mgt1} does not depend on $h$, which for the proposed setup should lead to $h-$independent convergence. Thus, while the convergence bounds can be shown to be independent of $h$, the practical convergence speed may still depend on $h$.

\subsubsection{\texorpdfstring{$h-$}{h-}independence}
We report the results for various $k$ using $h = 2^{-p}$, with $p = 5$ up to $p = 9$. 
\begin{table}[!hbt]
\centering
\begin{tabular}{cccc|ccc}
\hline
                              & \multicolumn{3}{c}{$k = 15$}      & \multicolumn{3}{c}{$k = 30$}     \\ \hline
\multicolumn{1}{c|}{$h$}      & $\nu = 1$ & $\nu = 2$ & $\nu = 4$ & $\nu = 1$ & $\nu = 2$ & $\nu = 4$ \\ \hline
\multicolumn{1}{c|}{$2^{-5}$} & 45        & 24        & 18        & $\odot$   & $\odot$   & $\odot$   \\
\multicolumn{1}{c|}{$2^{-6}$} & 34        & 22        & 18        & 66        & 37        & 28        \\
\multicolumn{1}{c|}{$2^{-7}$} & 36        & 22        & 18        & 52        & 33        & 27        \\
\multicolumn{1}{c|}{$2^{-8}$} & 40        & 24        & 18        & 54        & 34        & 27        \\
\multicolumn{1}{c|}{$2^{-9}$} & 42        & 23        & 18        & 58        & 36        & 27        \\ \hline
\end{tabular}
\label{tabhind}
\caption{Number of V-cycles for $k = 15$ and $k = 30$. $\nu$ denotes the number of {$\omega$-Jacobi} smoothing steps using $\omega = 4.5$ $\odot$ denotes the case where $kh > 0.625$, which is excluded.}
\end{table}
We observe that the solver indeed performs independently of $h$. {\color{black}In fact, it appears that with 4 smoothing steps, the solver converges in about $\tilde{C}k$ iterations, where $\tilde{C} \approx 1$ is a constant. The latter holds irrespective of the problem size.} 

\subsubsection{General results}
Now that we have established, both theoretically and numerically, that the convergence is $h-$independent once we take 4 smoothing steps, we revert back to letting $kh = 0.625$. Note that in order to obtain improved accuracy, $kh$ can be decreased without affecting the convergence behavior due to the $h-$independence. In \cref{tab:cons1} we report the results for MP 2-A using a constant wavenumber and Sommerfeld boundary conditions. 
\begin{table}[!hbt]
\centering
\caption{Number of V- $(\gamma =1)$ and W-cycles $(\gamma = 2)$ for constant $k$ (MP 2-A) using tol. $10^{-5}$. $\nu$ denotes the number of  {$\omega$-Jacobi} smoothing steps. $N_{D}$ is the size of the coarsest system.}
\scalebox{1}{
\begin{tabular}{c|cc|cc|cc|cc|cc}
& \multicolumn{2}{c|}{$k=50$} &  \multicolumn{2}{c|}{$k=100$} & \multicolumn{2}{c|}{$k=150$} &  \multicolumn{2}{c|}{$k=200$}  &  \multicolumn{2}{c}{$k=250$} \\ \hline
& \multicolumn{2}{c|}{$N = 6724$} &  \multicolumn{2}{c|}{$N = 26244$} & \multicolumn{2}{c|}{$N = 57600$} &  \multicolumn{2}{c|}{$N = 102400$} &  \multicolumn{2}{c}{$N = 160000$} \\ 
& \multicolumn{2}{c|}{$N_{D} = 8$} &  \multicolumn{2}{c|}{$N_{D} = 8$} & \multicolumn{2}{c|}{$N_{D} = 4$} &  \multicolumn{2}{c|}{$N_{D} =8$} &  \multicolumn{2}{c}{$N_{D} = 4$} \\ \hline
$\gamma$ & $1$   & $2$ & $1$   & $2$   & $1$    & $2$   & $1$    & $2$   & $1$   & $2$    \\ \hline 
$\nu =4$ & $58$  &$58$ & $104$ & $108$ & $155$  & $159$ & $209$  & $213$ & $267$ & $271$   \\
$\nu =5$ & $58$  &$58$ & $104$ & $104$ & $150$  & $166$ & $194$  & $229$ & $238$ & $287$   \\
$\nu =6$ & $55$  &$58$ & $99$  & $102$ & $139$  & $167$ & $183$  & $222$ & $226$ & $283$   \\
$\nu =7$ & $53$  &$60$ & $97$  & $101$ & $136$  & $163$ & $179$  & $219$ & $221$ & $280$   \\
$\nu =8$ & $53$  &$60$ & $95$  & $104$ & $131$  & $161$ & $178$  & $212$ & $218$ & $277$   \\
\end{tabular}}
\label{tab:cons1}
\end{table}
{\color{black}From \cref{tab:cons1} we observe that for the V-cycle the convergence improves if we use more smoothing steps.} Note that these additional smoothing steps are computationally cheap, as we use the $\omega-$Jacobi smoother. 
In general, we observe that we again need approximately $\tilde{C}k$ iterations with $\tilde{C} \approx 1$, in order to reach convergence. While the convergence behavior is promising, the results do provide some further insights into the behavior of the smoother. For example, we observe that moving from a V-cycle to a W-cycle does not seem to improve the performance. In fact, we observe that we need more iterations instead of less. One potential explanation could be that we know that in general the $\omega-$Jacobi smoother diverges for Helmholtz-type of problems. While this in itself does not have to lead to divergence (the smoothing property is still satisfied for a divergent smoother, see \cite{hackbusch2013multi}), the use of the W-cycle requires more smoothing iterations than the V-cycle. At some point we expect the efficiency of these smoothing steps to reduce, which could be reflected in the higher number of iterations. 

    \begin{figure*}[!hbt]
        \centering
        \begin{subfigure}
            \centering
            \includegraphics[width=0.49\textwidth]{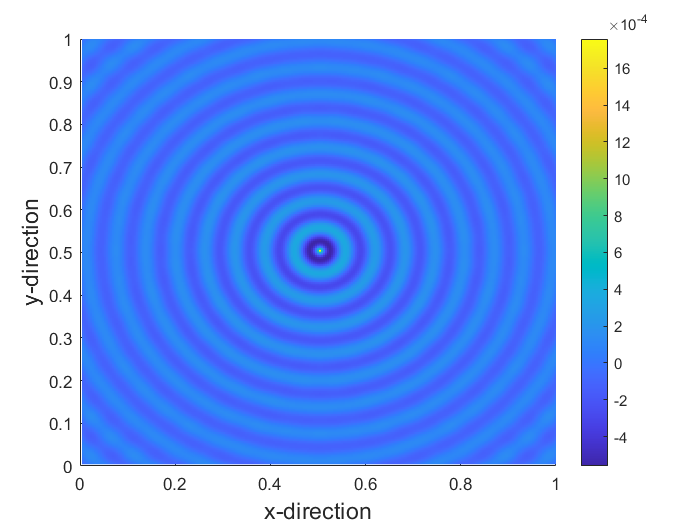}
            %\caption[Network2]%
            %{{\small $u(x,y)$ with $k(x,y)$ \textbf{sorted}}}    
            %\label{fignc75a}
        \end{subfigure}
        \hfill
        \begin{subfigure}
            \centering 
            \includegraphics[width=0.49\textwidth]{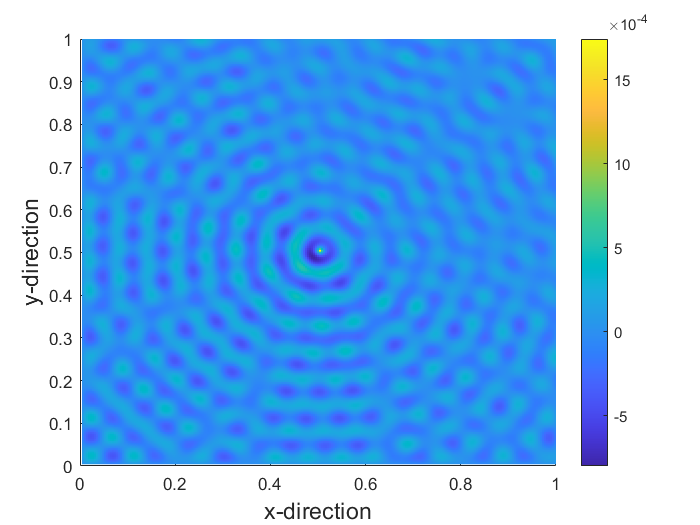}
            %\caption[]%
            %{{\small $k(x,y)$}}    
            %\label{fignc75b}
        \end{subfigure}
        \caption{Real part of the 2D solution for MP2-A $k = 100$ (left) and $k = 150$ (right)}
    \label{figab100}
    \end{figure*}

\subsection{2D results non-constant \texorpdfstring{$k(x,y)$}{k(x,y)}} \label{ncjacb}
In this subsection we report on the numerical results for the non-constant wavenumber model problem. We distinguish two cases; a medium-varying and a highly-varying wavenumber profile problem. In this respect the profile contains either smooth or sharp changes between a range of wavenumbers.
For an illustration of the wavenumber profile, see \cref{newnc1} for the smooth changing wavenumber profile and \cref{newnc2} for the profile where the wavenumber changes sharply. 
In all cases we make sure that the largest value of $k(x,y)$ still obeys the rule of thumb $kh \approx 0.625$. 

\FloatBarrier
\subsubsection{Smooth Changes}
\begin{figure*}[!hbt]
        \centering
        \begin{subfigure}
            \centering
            \includegraphics[width=0.49\textwidth]{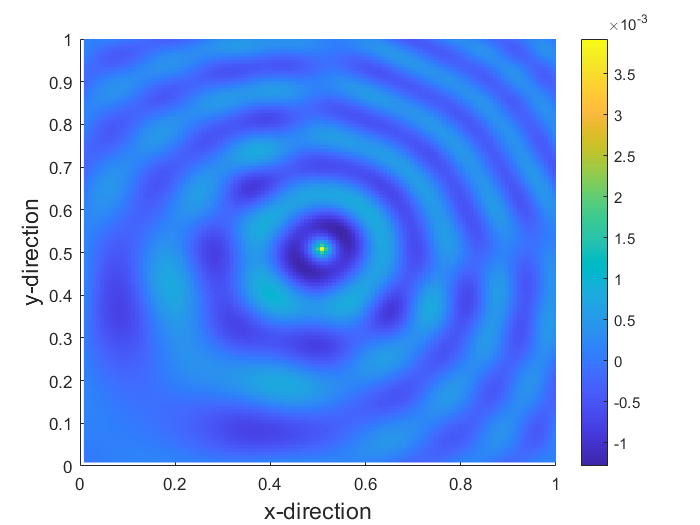}
            %\caption[Network2]%
            %{{\small $u(x,y)$ with $k(x,y)$ \textbf{sorted}}}    
            %\label{fignc75a}
        \end{subfigure}
        \hfill
        \begin{subfigure}
            \centering 
            \includegraphics[width=0.49\textwidth]{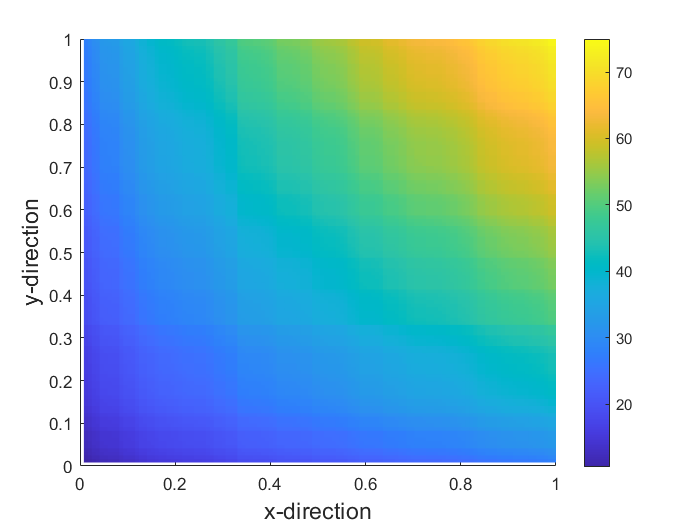}
            %\caption[]%
            %{{\small $k(x,y)$}}    
            %\label{fignc75b}
        \end{subfigure}
        \caption{Left: real part of the 2D numerical solution. Right: $k(x,y)$}
    \label{newnc1}
    \end{figure*}

\cref{newnc1} illustrates how the wavenumber varies between 10 and 75 for the medium-varying problem. Note that the transition from a low to high wavenumber goes gradually.
    \FloatBarrier
\cref{ncons1} contains the number of V-cycles needed to reach convergence for $k(x,y)$ varying between $(10,50)$ and $(10,75)$ respectively. We again observe that the number of iterations scales linearly with $k$. The iteration count follows approximately $\tilde{C}k$, with $\tilde{C} \approx 1$, where $k$ in this case denotes the largest wavenumber. In this event, the W-cycle does lead to a smaller number of iterations, compared to the case where we have a constant wavenumber. %This effect is most visible for $\nu = 4$ up to $\nu = 6$, in case $k(x,y)$ is between 10 and 50. Thus, as the wavenumber grows, the V-cycle is preferred over the W-cycle for this choice of smoothing scheme. 
\begin{table}[!h]
\centering
\captionsetup{width=.8 \textwidth}
\caption{Number of V- $(\gamma =1)$ and W-cycles $(\gamma = 2)$ for MP 2-B (medium variation). $\nu$ denotes the number of {$\omega$-Jacobi} smoothing steps. }
\scalebox{1}{
\begin{tabular}{ccc|cc}
\hline
                              & \multicolumn{2}{c}{$(k_1,k_2) = (10,50)$} & \multicolumn{2}{c}{$(k_1,k_2) = (10,75)$} \\ \hline
\multicolumn{1}{c|}{$\gamma$} & $1$                 & $2$                 & $1$                  & $2$                 \\ \hline
\multicolumn{1}{c|}{$\nu =4$} & $65$                & $60$                & $90$                & $88$               \\
\multicolumn{1}{c|}{$\nu =5$} & $62$                & $59$                & $86$                & $86$               \\
\multicolumn{1}{c|}{$\nu =6$} & $61$                & $58$                & $85$                & $85$               \\
\multicolumn{1}{c|}{$\nu =7$} & $60$                & $57$                & $84$                & $84$               \\
\multicolumn{1}{c|}{$\nu =8$} & $59$                & $57$                & $83$                & $83$               \\ \hline
\end{tabular}}
\label{ncons1}
\end{table}

\subsubsection{Sharp Changes}
\begin{figure*}[!hbt]
        \centering
        \begin{subfigure}
            \centering
            \includegraphics[width=0.49\textwidth]{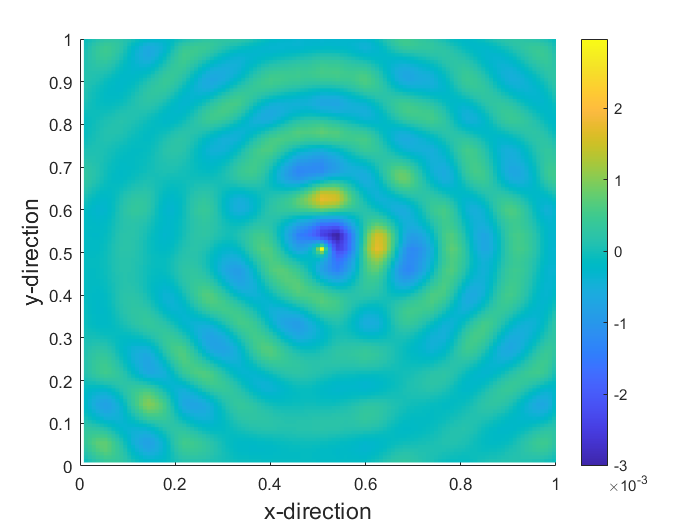}
            %\caption[Network2]%
            %{{\small $u(x,y)$ with $k(x,y)$ \textbf{sorted}}}    
            %\label{fignc75a}
        \end{subfigure}
        \hfill
        \begin{subfigure}
            \centering 
            \includegraphics[width=0.49\textwidth]{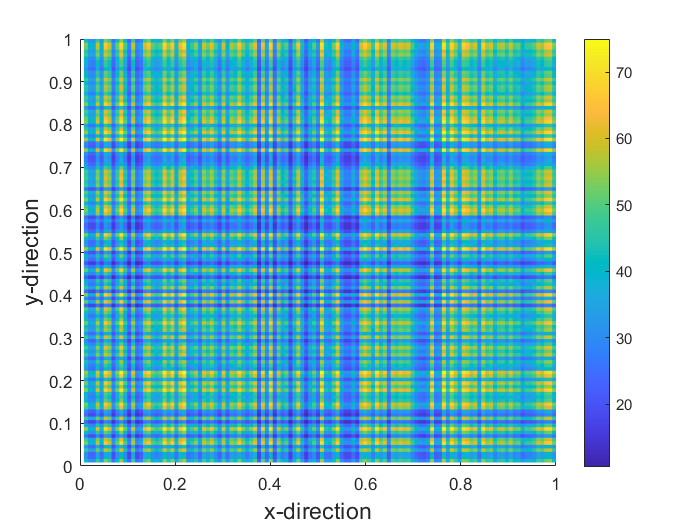}
            %\caption[]%
            %{{\small $k(x,y)$}}    
            %\label{fignc75b}
        \end{subfigure}
        \caption{Left: real part of the 2D numerical solution. Right: $k(x,y)$}
    \label{newnc2}
    \end{figure*}

In this subsection we report on the results for the highly-varying profile model problem, containing sharp changes between the wavenumbers $k=10$ and $k=75$. \cref{newnc2} illustrates how the wavenumber varies between 10 and 75 for the medium-varying problem. In this case, we let $k(x,y)$ vary randomly, which gives a wavenumber which varies highly across the entire numerical domain. 
    \FloatBarrier
\cref{ncons2} contains the number of V-cycles needed to reach convergence for $k(x,y)$ varying between $(10,50)$ and $(10,75)$ respectively. We immediately observe that the number of iterations goes up by a factor of approximately 1.5 once we allow for the sharp and rapid changes in the wavenumber. For example for the medium-varying case containing more smooth transition between the wavenumbers, using 8 smoothing steps, we had 59 and 57 iterations respectively for the V- and W-cycle. In this case, we have 94 iterations for both the V- and W-cycle. Another interesting observation is that for the W-cycle, $k(x,y)$ varying between 10 and 75 leads to a larger number of iterations when $\nu > 4$. This could again be an indication of the smoother losing some of its efficacy when we allow for more smoothing operations using the W-cycle. 
Also for this model problem, we conclude that as the wavenumber grows, the V-cycle is preferred over the W-cycle for this choice of smoothing scheme. In the reported cases, the method converges in approximately $~1.5k$ iterations, where $k$ is the largest admissible wavenumber. 
\begin{table}[!h]
\centering
\scalebox{1}{
\begin{tabular}{ccc|cc}
\hline
                              & \multicolumn{2}{c}{$(k_1,k_2) = (10,50)$} & \multicolumn{2}{c}{$(k_1,k_2) = (10,75)$} \\ \hline
\multicolumn{1}{c|}{$\gamma$} & $1$       & \multicolumn{1}{c|}{$2$}      & $1$                 & $2$                 \\ \hline
\multicolumn{1}{c|}{$\nu =4$} & $102$     & \multicolumn{1}{c|}{$96$}    & $111$               & $107$               \\
\multicolumn{1}{c|}{$\nu =5$} & $97$     & \multicolumn{1}{c|}{$95$}    & $103$               & $105$               \\
\multicolumn{1}{c|}{$\nu =6$} & $95$     & \multicolumn{1}{c|}{$95$}    & $101$               & $104$               \\
\multicolumn{1}{c|}{$\nu =7$} & $94$     & \multicolumn{1}{c|}{$94$}    & $102$               & $104$               \\
\multicolumn{1}{c|}{$\nu =8$} & $94$     & \multicolumn{1}{c|}{$94$}    & $102$               & $104$               \\ \hline
\end{tabular}}
\captionsetup{width=.8 \textwidth}
\caption{Number of V- $(\gamma =1)$ and W-cycles $(\gamma = 2)$ for MP 2-B (high variation). $\nu$ denotes the number of {$\omega$-Jacobi} smoothing steps.}
\label{ncons2}
\end{table}

\subsection{GMRES-smoothing} In this section we explore how the solver performs when we use GMRES(3) smoothing. In \cite{cools2014new}, where the authors propose coarsening on the CSL instead of the original Helmholtz operator, GMRES(3) is used as the only smoother. The motivation for this lies in the fact that indeed the $\omega-$Jacobi diverges when combined with the original and standard configuration of the multigrid method.
By using a few GMRES iterations, the smoother acts as a polynomial smoother. Already in \cite{elman2001multigrid}, the use of GMRES(3) was proposed as a smoother. There, the authors manually optimized the smoothing schedule, where the number of smoothing steps ranges between 13 and 39 smoothing steps for a constant $k$ ranging from $4\pi$ to $32\pi$. More generally in the literature, on each level where the Jacobi smoother becomes unstable, the smoother is replaced by the Krylov iterations ranging between 5-40 iterations per smoothing step \cite{elman2001multigrid,chen2014robust, chen2015multilevel}. In this work we let the number of smoothing steps using GMRES(3) vary between 1 and 5. 
One advantage of this approach is that by allowing for only 3 GMRES iterations per smoothing step, the cost of applying GMRES does not increase significantly with respect to memory and computation time. We perform the same numerical experiments using the model problems mentioned previously. 

\subsubsection{2D results constant \texorpdfstring{$k$}{k}}
In this subsection we report the results for the constant wavenumber model problem (MP 2-A). As we are coarsening using the CSL, we want to distinguish between the cases where the complex shift is large versus small. We therefore start by keeping the shift $\beta_2 = 0.7$, similar to the one used in all previous experiments. Results are reported in \cref{tab:cons1g} for 1 up to 5 smoothing steps. We immediately observe that the number of iterations for both the V- and W-cycle are drastically reduced. If we compare these results to the ones obtained when using the $\omega-$Jacobi smoother in \cref{tab:cons1}, the reduction in the number of iterations is approximately 2.5 times. For example, when using 5 smoothing steps within a V-cycle, the $\omega-$Jacobi smoother requires 238 iterations to reach convergence. When GMRES(3) is used as a smoother with a similar amount of smoothing steps, 88 iterations are required to reach convergence. Note that there appears to be no difference in the number of iterations when using the V-cycle and W-cycle. 
\begin{table}[!hbt]
\centering
\caption{Number of V- $(\gamma =1)$ and W-cycles $(\gamma = 2)$ for constant $k$ (MP 2-A) using tol. $10^{-5}$. $\nu$ denotes the number of {GMRES(3)} smoothing steps with $\beta_2 = 0.7$}
\scalebox{1}{
\begin{tabular}{c|cc|cc|cc|cc|cc}
& \multicolumn{2}{c|}{$k=50$} &  \multicolumn{2}{c|}{$k=100$} & \multicolumn{2}{c|}{$k=150$} &  \multicolumn{2}{c|}{$k=200$}  &  \multicolumn{2}{c}{$k=250$} \\ \hline
& \multicolumn{2}{c|}{$N = 6724$} &  \multicolumn{2}{c|}{$N = 26244$} & \multicolumn{2}{c|}{$N = 57600$} &  \multicolumn{2}{c|}{$N = 102400$} &  \multicolumn{2}{c}{$N = 160000$} \\ 
& \multicolumn{2}{c|}{$N_{D} = 8$} &  \multicolumn{2}{c|}{$N_{D} = 8$} & \multicolumn{2}{c|}{$N_{D} = 4$} &  \multicolumn{2}{c|}{$N_{D} =8$} &  \multicolumn{2}{c}{$N_{D} = 4$} \\ \hline
$\gamma$ & $1$   & $2$ & $1$    & $2$   & $1$    & $2$   & $1$    & $2$   & $1$   & $2$    \\ \hline 
$\nu =1$ & $37$  &$36$ & $68$	& $67$ & $99$  & $98$ & $132$  & $131$ & $162$ & $161$   \\
$\nu =2$ & $29$  &$29$ & $53$	& $53$ & $78$  & $78$ & $104$  & $104$ & $128$ & $128$   \\
$\nu =3$ & $24$  &$24$ & $45$  & $45$ & $67$  & $67$ & $89$  & $89$ & $112$ & $112$   \\
$\nu =4$ & $22$  &$22$ & $40$  & $40$ & $59$  & $59$ & $78$  & $78$ & $98$ & $98$   \\
$\nu =5$ & $20$  &$20$ & $36$  & $36$ & $53$  & $53$ & $71$  & $71$ & $88$ & $88$   \\
\end{tabular}}
\label{tab:cons1g}
\end{table}
\FloatBarrier
Next, we repeat the same numerical experiments, however this time we set the complex shift to $\beta_2 = k^{-1}$. Note that this is a very small shift, in fact the CSL matrix starts being a close resemblance of the original Helmholtz operator. Results are reported in \cref{tab:cons2g}. For all cases, we report that there is a drastic improvement in the number of iterations. {\color{black}We obtain a solver which is close to $k-$ independent when we use $\nu >2$ and $\gamma = 2$}. Unlike the previous case where $\beta_2 = 0.7$ in \cref{tab:cons1g}, we now do observe a lower number of iterations for the W-cycle.

In fact, additional numerical experiments confirm that the results are similar if we let $\beta_2 = 0$, which results in the original Helmholtz operator. Thus, when using GMRES(3) as a smoother, the original Helmholtz operator can in fact be used instead of the CSL for coarsening within the multigrid hierarchy. However, the fast convergence is only observed when combined with the high-order prolongation and restriction operators. 

Given that we obtained the best numerical results with $\beta_2 = k^{-1}$ without any additional costs, we continue with this shift for the CSL in the upcoming sections. 
\begin{table}[!hbt]
\centering
\scalebox{1}{
\begin{tabular}{c|cc|cc|cc|cc|cc}
& \multicolumn{2}{c|}{$k=50$} &  \multicolumn{2}{c|}{$k=100$} & \multicolumn{2}{c|}{$k=150$} &  \multicolumn{2}{c|}{$k=200$}  &  \multicolumn{2}{c}{$k=250$} \\ \hline
& \multicolumn{2}{c|}{$N = 6724$} &  \multicolumn{2}{c|}{$N = 26244$} & \multicolumn{2}{c|}{$N = 57600$} &  \multicolumn{2}{c|}{$N = 102400$} &  \multicolumn{2}{c}{$N = 160000$} \\ 
& \multicolumn{2}{c|}{$N_{D} = 8$} &  \multicolumn{2}{c|}{$N_{D} = 8$} & \multicolumn{2}{c|}{$N_{D} = 4$} &  \multicolumn{2}{c|}{$N_{D} =8$} &  \multicolumn{2}{c}{$N_{D} = 4$} \\ \hline
$\gamma$ & $1$   & $2$ & $1$    & $2$   & $1$    & $2$   & $1$    & $2$   & $1$   & $2$    \\ \hline 
$\nu =1$ & $14$  &$7$ & $24$	& $10$ & $39$  & $19$ & $51$  & $24$ & $64$ & $29$   \\
$\nu =2$ & $8$  &$5$ & $13$	& $7$ & $22$  & $10$ & $28$  & $13$ & $34$ & $16$   \\
$\nu =3$ & $6$  &$5$ & $10$  & $6$ & $16$  & $9$ & $20$  & $10$ & $24$ & $12$   \\
$\nu =4$ & $6$  &$5$ & $8$  & $5$ & $12$  & $7$ & $15$  & $9$ & $18$ & $10$   \\
$\nu =5$ & $5$  &$5$ & $7$  & $5$ & $11$  & $7$ & $13$  & $8$ & $15$ & $9$   \\
\end{tabular}}
\caption{Number of V- $(\gamma =1)$ and W-cycles $(\gamma = 2)$ for constant $k$ (MP 2-A) using tol. $10^{-5}$. $\nu$ denotes the number of {GMRES(3)} smoothing steps with $\beta_2 = k^{-1}$.}
\label{tab:cons2g}
\end{table}
\FloatBarrier
%\vspace{-0.5cm}
\subsubsection{2D results non-constant \texorpdfstring{$k(x,y)$}{k(x,y)}}
Finally, we investigate the convergence behavior using the GMRES(3) smoother on the highly-varying problem (MP 2-B2). In \cref{ncjacb}, we observed that this was the hardest problem to solve in terms of heterogeneous problem, due to the wavenumber varying highly between $k = 10$ and $k = 75$ across the numerical domain. Thus, we only test for this case in this subsection. 
\FloatBarrier
\subsubsection{Sharp Changes}
Results are reported in \cref{tab:ncons2g} and again indicate a drastic improvement compared to the case where we use the $\omega-$Jacobi smoother. In fact, even for this highly-varying model containing sharp changes between the wavenumbers $k = 10$ and $k=75$, where the wavenumber is allowed to vary randomly between 10 and 75 across the domain, we obtain a $k-$ and $h-$independent multigrid solver for $\nu \geq 3$. Similar to the previous case, even when using the original Helmholtz operator, we reach convergence when using GMRES(3) as a smoother. However, the $k-$independence is only observed when the higher-order prolongation and restriction scheme is used. 
{\color{black}These results are promising and provide a solid framework for future research, given that the current industry standard (CSL preconditioner with multigrid inversion) works well for homogeneous problems, but is less suitable for heterogeneous problems with sharp changes between wavenumbers due to the instability of the inexact inversion using geometric multigrid \cite{sheikh2016accelerating,erlangga2008multilevel}}. 
\begin{table}[!hbt]
\centering
\scalebox{1}{
\begin{tabular}{ccc|cc}
\hline
                              & \multicolumn{2}{c}{$(k_1,k_2) = (10,50)$} & \multicolumn{2}{c}{$(k_1,k_2) = (10,75)$} \\ \hline
\multicolumn{1}{c|}{$\gamma$} & $1$       & \multicolumn{1}{c|}{$2$}      & $1$                 & $2$                 \\ \hline
\multicolumn{1}{c|}{$\nu =1$} & $28$     & \multicolumn{1}{c|}{$12$}    & $31$               & $12$               \\
\multicolumn{1}{c|}{$\nu =2$} & $16$     & \multicolumn{1}{c|}{$8$}    & $17$               & $7$               \\
\multicolumn{1}{c|}{$\nu =3$} & $12$     & \multicolumn{1}{c|}{$7$}    & $12$               & $6$               \\
\multicolumn{1}{c|}{$\nu =4$} & $10$     & \multicolumn{1}{c|}{$6$}    & $10$               & $6$               \\
\multicolumn{1}{c|}{$\nu =5$} & $9$     & \multicolumn{1}{c|}{$6$}    & $9$               & $6$               \\ \hline
\end{tabular}}
\captionsetup{width=.8 \textwidth}
\caption{Number of V- $(\gamma =1)$ and W-cycles $(\gamma = 2)$ for MP 2-B (high variation). $\nu$ denotes the number of {GMRES(3)} smoothing steps and $\beta_2 = {k_{\max}}^{-1}$.}
\label{tab:ncons2g}
\end{table}
\FloatBarrier

\section{Conclusion} \label{sec6}
In this paper, we present the \textit{first} full multigrid hierarchy for both the V- and W-cycles using classic multigrid components for the indefinite nonsymmetric Helmholtz problem. Here, we emphasized the use of classical components, such as the weighted Jacobi smoother. The resulting algorithm shows $h-$independent convergence and thus adheres to the classic multigrid features. Two novel and striking features should be mentioned. First of all, no restriction is imposed on the number of grid points on the coarsest grid. Second of all, no level-dependent parameters are needed. 

Apart from the numerical results, we provide a new theory to assess convergence of highly indefinite linear systems. Where most proofs require some adjointness and/or symmetry assumptions on the underlying linear system, we have constructed theoretical notions which do not require such assumptions. In fact, we have found that the addition of a complex shift, solely for the purpose of coarsening in the multigrid hierarchy, combined with higher-order interpolation schemes for the inter-grid transfer operators leads to an HPD system. The positive definiteness can be verified for medium-sized wavenumbers by assessing properties of the matrix. {\color{black}As a result, it can be shown that the convergence is independent of $h$}. Thus, the computation of the spectral radius can be circumvented, especially since the analytical eigenvalues can only be computed in case of Dirichlet boundary conditions. 
We also study the behavior of the solver once we use GMRES(3) as a smoother instead of $\omega-$Jacobi. The convergence significantly improves and less smoothing steps are needed. Second of all, if we use GMRES(3) as a smoother, we can keep the complex shift very small without paying a penalty in terms of additional costs. As a result, we obtain close to $k-$ and $h-$independent convergence, both for the constant wavenumber model problem and the non-constant wavenumber model problem. In fact, we can even use the original Helmholtz operator for coarsening in this configuration without losing efficiency and effectiveness as the GMRES iterations provide polynomial smoothing. 

This work is one of the few to demonstrate, both theoretically and numerically, the existence and potential of multigrid methods for such linear systems. However, a lot of future work remains to be done. For example, future research could focus on exploring different combinations of interpolation schemes and smoothers and extending proofs. Also, instead of coarsening on the CSL, other candidates in terms of 'near positive definite operators' could also be explored. We have shown convergence for $\omega-$Jacobi, mainly because of its simplicity and parallelizability, but many other options remain to be examined. Moreover, once a convergent solver has been established, the next step in shifting the paradigm would be to study wavenumber independent convergence. 

\bibliographystyle{siamplain}
%\bibliography{references3}

%%% BIB %%%

%%%%%%%%%%%

\end{document}

%% file: ex_shared.tex
\usepackage{lipsum}
\usepackage{amsfonts}
\usepackage{graphicx}
\usepackage{epstopdf}
\usepackage{algorithmic}
\ifpdf
  \DeclareGraphicsExtensions{.eps,.pdf,.png,.jpg}
\else
  \DeclareGraphicsExtensions{.eps}
\fi

\newsiamremark{remark}{Remark}
\newsiamremark{hypothesis}{Hypothesis}
\crefname{hypothesis}{Hypothesis}{Hypotheses}
\newsiamthm{claim}{Claim}

\headers{Stand-alone Multigrid for Helmholtz}{Vandana Dwarka, Cornelis Vuik}

\title{Stand-alone multigrid for Helmholtz revisited: towards convergence using standard components}
\author{Vandana Dwarka\thanks{Department of Applied Mathematics, Delft University of Technology, Delft, the Netherlands
  (\email{v.n.s.r.dwarka@tudelft.nl}).}
\and Cornelis Vuik \thanks{Department of Applied Mathematics, Delft University of Technology, Delft, the Netherlands
  (\email{c.vuik@tudelft.nl}, \url{homepage.tudelft.nl/d2b4e/}).}}

\usepackage{amsopn}